\documentclass[11pt]{amsart}
\usepackage{amsfonts,amssymb}
\usepackage{eucal}
\usepackage{graphics}
\usepackage{color}

\newtheorem{theorem}{Theorem}
\newtheorem{lemma}[theorem]{Lemma}
\newtheorem{prop}[theorem]{Proposition}

\newtheorem{corollary}[theorem]{Corollary}

\newcommand{\C}{\mathbb C}
\newcommand{\R}{\mathbb R}
\renewcommand{\P}{\mathbb P}
\renewcommand{\O}{\EuScript O}

\newcommand{\D}{\mathbb D}
\newcommand{\F}{\mathbb F}
\newcommand{\A}{\mathbb A}
\newcommand{\Z}{\mathbb Z}

\newcommand{\iso}{\cong}

\renewcommand{\o}{\omega}
\newcommand{\suchthat}{\;:\;}
\newcommand{\wigglygeq}{\gtrsim}
\newcommand{\tildei}{\tilde{\text{\it\i}}}

\title[Ramanujam's surface]{The symplectic topology
of\\ Ramanujam's surface}
\author{Paul Seidel, Ivan Smith}
\date{13 November 2004. Keywords: Stein manifold, contractible affine
  surface, exotic symplectic structure, subcritical Stein manifold.
  }

\begin{document}

\begin{abstract} Ramanujam's surface $M$ is a contractible affine
  algebraic surface which is not homeomorphic to the affine plane.
  For any $m>1$ the product $M^m$ is diffeomorphic to Euclidean space
  $\R^{4m}$. We show that, for every $m>0$, $M^m$ cannot be
  symplectically embedded into a subcritical Stein manifold.  This
  gives the first examples of exotic symplectic structures on
  Euclidean space which are convex at infinity.  It follows that any
  exhausting plurisubharmonic Morse function on $M^m$ has at least
  three critical points, answering a question of Eliashberg.  The
  heart of the argument involves showing a particular Lagrangian torus
  $L$ inside $M$ cannot be displaced from itself by any Hamiltonian
  isotopy, via a careful study of pseudoholomorphic discs with
  boundary on $L$.
\end{abstract}

\maketitle

\section{Introduction}

Ramanujam showed in \cite{ramanujam71} that the complement $M$ of a
certain singular curve in the Hirzebruch surface $\F_1$ is a
\emph{contractible} algebraic surface.  Algebro-geometrically, $M$
is distinguished from the affine plane $\A^2$ by being of log
general type (having log Kodaira dimension $2$, cf.
\cite{miyanishi}). Topologically, in spite of being contractible,
$M$ is not homeomorphic to $\R^4$, since its fundamental group at
infinity is nontrivial. Now consider the $m$-fold product $M^m = M
\times \cdots \times M$. This is still of log general type, in
particular not isomorphic to $\A^{2m}$ as an algebraic variety.
However, for $m \geq 2$ the fundamental group at infinity becomes
trivial, and as a consequence $M^m$ is diffeomorphic to $\R^{4m}$;
indeed, Dimca [unpublished] observed that any contractible affine
variety of complex dimension $d \geq 3$ is diffeomorphic to
$\R^{2d}$ (cf. \cite[Theorem 3.2]{zaidenberg98}). The conclusion is
that the smooth manifolds $\R^{4m}$, $m \geq 2$, admit nonstandard
algebraic variety structures.

The aim of the present note is to consider this phenomenon from the
symplectic perspective. We will equip $M$ with an exhausting
plurisubharmonic function $\phi_M$, which makes it into a Stein
manifold of finite type, and consider the associated symplectic form
$\o_M$ (for this and related terminology, see Sections \ref{sec:def}
and \ref{sec:essential}). On $M^m$ we take the product structure.

\begin{theorem} \label{th:main}
For all $m$, $M^m$ cannot be symplectically embedded into a
subcritical Stein manifold.
\end{theorem}

In particular, for $m \geq 2$ we obtain a symplectic structure on
$\R^{4m}$ which is exotic (in the usual sense, of not admitting an
embedding into the standard $\R^{4m}$). Any sufficiently large
relatively compact part of $M$ is exotic in the same sense, because
of the finite type property (cf. Lemma \ref{th:weinstein-embed} and
Lemma \ref{th:stupid} below). There is a consequence which can be
stated purely in terms of Stein geometry, answering a question of
Eliashberg \cite[Problem 3]{eliashberg94}:

\begin{corollary} \label{th:cor}
For all $m$, any exhausting plurisubharmonic Morse function on $M^m$
must have at least $3$ critical points.
\end{corollary}

It seems appropriate to compare Theorem \ref{th:main} with some
other known results. There are several constructions of exotic
symplectic structures on $\R^{2n}$ for $n \geq 2$, starting with the
abstract existence theorem of \cite[Corollary 0.4.$A_2'$]{gromov85}.
However, in contrast to our example, the resulting symplectic forms
are not known to be convex at infinity; in fact, at least one
construction \cite{muller90} is explicitly designed to violate that
condition.   In a somewhat different direction, we should mention
that Eliashberg \cite{eliashberg94} has given candidates, by an
explicit Lagrangian handle decomposition, for Stein
\emph{subdomains} of $\C^{2n}$, $n>2$, which are diffeomorphic to
balls, and for which he conjectures that the conclusion of Corollary
\ref{th:cor} still holds.

The main result of \cite{ramanujam71} asserts that $\A^2$ is the
only algebraic surface which is contractible {\em and}
simply-connected at infinity. The symplectic counterpart of this is
the observation that (assuming the $3$-dimensional Poincar{\'e}
conjecture) any $4$-dimensional Weinstein manifold which is
contractible, simply connected at infinity, complete, and of finite
type, is symplectically isomorphic to standard $\R^4$.  The proof
relies on the uniqueness of tight contact structures on $S^3$
\cite{eliashberg92} and the description of Stein fillings of this
structure via families of holomorphic discs \cite{eliashberg90}.
Moreover, the picture changes if one drops the finite type
condition: Gompf \cite{gompf98} has used a suitable infinite
handlebody decomposition to produce Stein structures on uncountably
many manifolds homeomorphic, but not diffeomorphic, to $\R^4$.

The essential ingredient in our proof of Theorem \ref{th:main} is a
particular Lagrangian torus $L \subset M$, described below. By a
careful study of pseudo-holomorphic discs, and invoking a theorem of
Chekanov \cite{chekanov95}, we show that $L \subset M$ cannot be
displaced from itself by a Hamiltonian isotopy. More generally, if
$i: M \rightarrow N$ is a symplectic embedding into a complete Stein
manifold, then $i(L) \subset N$ has the same non-displaceability
property. This, together with the corresponding facts for products
$L^m$, leads easily to Theorem \ref{th:main}.

The heart of the argument involves considering two different
compactifications $X$ and $\bar{X}\cong\F_1$ of $M$.  The complement
$\bar{X}\backslash M$ is a curve $\bar{S}=\bar{S}'\cup\bar{S}''$
with two irreducible components, one of which $\bar{S}''$ has a cusp
singularity; the Lagrangian torus $L$ lies in a neighbourhood of the
cusp.  Rather than working with a single Lagrangian torus, we
consider a family $\{L_t\}$ which as $t\rightarrow 0$ collapses into
the cusp point on $\bar{S}''$.  Intuitively, the limit of a family
of holomorphic discs $\{(\D,\partial \D)\rightarrow (M,L_t)\}_{t\in
(0,1]}$ as $t\rightarrow 0$ is either a holomorphic sphere in
$\bar{X}$ disjoint from $\bar{S}'$, or is the constant map to the
cusp point.  The first case is excluded since $\bar{S}'$ is ample;
the second is impossible for topological reasons concerning the
fundamental group $\pi_1(\bar{V}\backslash\bar{S}'')$, where
$\bar{V}\subset\bar{X}$ is a small neighbourhood of the cusp.  The
upshot is that no such families of discs can exist, enabling us to
appeal to Chekanov's work.

For technical reasons, we in fact work with a blown-up
compactification $X \rightarrow \bar{X}$ of $M$ in which the
complement $S = X \backslash M$ is a divisor with normal crossings.
The tori $L_t$ now appear as so-called linking tori for a normal
crossing point $p$ of $S$. To construct them as manifestly
Lagrangian tori, and to make the abovementioned limiting argument
for holomorphic discs rigorous, we use a simple algebro-geometric
trick. Take $\C P^1 \times X$ and blow up the point $(0,p)$,
obtaining a threefold $Y$ with a projection $Y \rightarrow \C P^1$.
The singular fibre $Y_0$ has an irreducible component which is a $\C
P^2$. We take a Clifford torus $K_0$ in that component, and move it
by parallel transport to obtain family of Lagrangian tori $K_t$, $t
\in [-1;1]$, in the nearby fibres $Y_t$. For $t \neq 0$, these
fibres are naturally identified with $X$, and we define $L_t$ to be
the image of $K_t$ under this identification, for $t \in (-1;0]$.
Since the total submanifold $K = \bigcup_{t \in [-1,1]}K_t \subset
Y$ is Lagrangian, Gromov compactness can be applied directly to
families of discs with boundary in $K$. The drawback of this
argument is that the $Y_t \iso X$ carry varying K{\"a}hler forms. To
take account of this, we give a careful discussion of Stein
deformations in Section~\ref{sec:def}, and introduce in
Section~\ref{sec:essential} the technical notion of a
``Stein-essential'' Lagrangian submanifold. The idea is that for any
given $E>0$, we can deform our Stein structure and our Lagrangian
submanifold, in such a way that at the endpoint of the deformation,
there are no pseudo-holomorphic discs of area less than $E$.
Section~\ref{sec:linking} introduces Lagrangian linking tori, and
describes the implications of a linking torus being Stein-essential.
Only in Section~\ref{sec:conclusion} are these ingredients assembled
to derive Theorem~\ref{th:main}.

There are at least two possible alternative ways of analyzing the
symplectic nature of $M^m$. One could try to use Floer homology or
symplectic homology as introduced by Viterbo
\cite{viterbo97a,viterbo97b} and Cieliebak-Floer-Hofer
\cite{cieliebak-floer-hofer95}.  In fact, as we intend to discuss
elsewhere, existence of a Stein-essential Lagrangian submanifold
already implies that $SH^*(M^m)\neq 0$. To compute $SH^*(M^m)$
precisely would presumably require an analysis of the Reeb flow at
infinity, though since Floer homology behaves well under products
\cite{oancea04} it would be enough to do the computation for $M$
itself. Floer homology may also distinguish between different exotic
symplectic structures on $\R^{4m}$, which falls outside the scope of
the arguments used here. The other possible approach would be via a
symplectic field theory decomposition argument. The aim would be to
prove that if $M^m$ is subcritical, there has to be a non-constant
algebraic map $\A^1 \rightarrow M$, contradicting a property of
contractible surfaces of log general type \cite[Theorem
4.7.1]{miyanishi}. Eliashberg has announced a theorem which says
that if a smooth projective variety contains a smooth ample divisor
with subcritical complement, then the variety has many rational
curves.  In examples, it appears that these rational curves are
closures of maps of $\A^1$ to the complement, but this is not
well-understood in general. Closely related results have been
obtained by Biolley \cite{ALB}.

{\em Acknowledgments.} The first author would like to thank Denis
Auroux for an illuminating discussion.  Thanks go to Kai Cieliebak
for helpful comments on an earlier version of the paper.  This
research was partially supported by NSF grant DMS-0405516 and a
grant from the Nuffield foundation NUF-NAL/00876/G.

\section{Background\label{sec:def}}

We begin by reviewing the definitions and some elementary results.
This follows \cite{eliashberg-gromov93, eliashberg94,
biran-cieliebak01} with some modifications. The proofs have been
relegated to the Appendix.

Take a manifold $M$ equipped with a symplectic form $\o_M$, a
one-form $\theta_M$ such that $d\theta_M = \o_M$, and an exhausting
(which means proper and bounded below) smooth function $\phi_M: M
\rightarrow \R$. Let $\lambda_M$ be the Liouville vector field
associated to $\theta_M$, so $\o_M(\lambda_M,\cdot) = \theta_M$. The
quadruple $(M,\o_M,\theta_M,\phi_M)$ is a {\em convex symplectic
manifold} if there is a sequence $c_1 < c_2 <\dots$ converging to
$+\infty$, such that $d\phi_M(\lambda_M) > 0$ on each level set
$\phi_M^{-1}(c_k)$. We call a convex symplectic manifold {\em
complete} if the flow of $\lambda_M$ exists for all positive times
(the corresponding statement for negative times is always true), and
{\em of finite type} if there is a $c_0$ such that
$d\phi_M(\lambda_M) > 0$ on $\phi_M^{-1}([c_0;+\infty))$. Note that
if $M$ is complete and of finite type, then the flow of $\lambda_M$
defines a diffeomorphism $f: [0;\infty) \times \phi_M^{-1}(c_0)
\rightarrow \phi_M^{-1}([c_0;+\infty))$ satisfying $f^*\theta_M =
e^r(\theta_M|\phi_M^{-1}(c_0))$, where $r$ is the variable in
$[0;\infty)$. Hence $M$ is a symplectic manifold with a conical end.

\begin{lemma} \label{th:weinstein-embed}
Let $M,N$ be convex symplectic manifolds, with $M$ of finite type
and $N$ complete. Take $c_0$ such that $d\phi_M(\lambda_M) > 0$ on
$\phi_M^{-1}([c_0;+\infty))$. Then any embedding $i:
\phi_M^{-1}((-\infty;c_0]) \rightarrow N$ such that $i^*\theta_N -
\theta_M$ is an exact one-form, can be extended to an embedding $M
\rightarrow N$ with the same property.
\end{lemma}

Let $(\o_{M,t},\theta_{M,t},\phi_{M,t})$, $0 \leq t \leq 1$, be a
smooth family of convex symplectic structures on a fixed manifold
$M$. We say that this is a {\em convex symplectic deformation} if
the following two additional conditions hold: the function $(t,x)
\mapsto \phi_{M,t}(x)$ on $[0;1] \times M$ is proper; and for each
$t \in [0;1]$ there is a neighbourhood $t \in I \subset [0;1]$ and a
sequence $c_1 < c_2 <\dots$ converging to $+\infty$, such that
$d\phi_{M,s}(\lambda_{M,s}) > 0$ along $\phi_{M,s}^{-1}(c_k)$ for
all $k$ and all $s \in I$. A convex symplectic deformation is called
{\em complete} if all the convex symplectic structures in it are
complete, and {\em of finite type} if there is a $c_0$ such that
$d\phi_{M,t}(\lambda_{M,t}) > 0$ on $\phi_{M,t}^{-1}([c_0;+\infty))$
for all $t$.

\begin{lemma} \label{th:weinstein-deform}
Let $(\o_{M,t},\theta_{M,t},\phi_{M,t})$ be a complete convex
symplectic deformation. For any relatively compact open subset $U
\subset M$, there is a smooth family of embeddings $j_t: U
\rightarrow M$ starting with $j_0 = id$, such that
$j_t^*\theta_{M,t} - \theta_{M,0}$ are exact one-forms on $U$.
\end{lemma}

\begin{lemma} \label{th:weinstein-deform-2}
Let $(\o_{M,t},\theta_{M,t},\phi_{M,t})$ be a complete finite type
convex symplectic deformation. Then there is a smooth family of
diffeomorphisms $f_t: M \rightarrow M$, starting with $f_0 = id$,
such that $f_t^*\theta_{M,t} - \theta_{M,0} = dR_t$, where $(t,x)
\mapsto R_t(x)$ is a compactly supported function on $[0;1] \times
M$.
\end{lemma}

It may be instructive to compare our definitions with some that
appear elsewhere in the literature.  The condition of a manifold
being \emph{Weinstein}, defined in \cite{eliashberg94}, is related
to but stronger than being convex symplectic; the function $\phi_M$
is tied closely to $\lambda_M$ by a Lyapunov condition, which is
somewhat more restrictive than our requirements. Convex symplectic
manifolds were introduced in \cite{eliashberg-gromov93}, defined as
exact symplectic manifolds $(M,\o_M,\theta_M)$ together with an
exhaustion by relatively compact subsets $U_1 \subset \bar{U}_1
\subset U_2 \subset \bar{U}_2 \subset U_3 \cdots$, such that each
$\partial U_k$ is a smooth hypersurface and convex of contact type.
This co-incides with our notion, but we choose to describe the
exhaustion $U_k = \phi_M^{-1}((-\infty;c_k))$ via sublevel sets of
some function.  Our definition of convex symplectic deformation
stays close to the same picture, since locally in the deformation
parameter $t$, the manifolds $(M,\o_{M,t},\theta_{M,t})$ have
smoothly varying exhaustions $U_{k,t} =
\phi_{M,t}^{-1}((-\infty;c_k))$.  For convex symplectic manifolds
which are both finite type and complete, this is the natural
analogue of the notion of deformation in \cite{eliashberg94}.

A {\em Stein manifold} $(M,J_M,\phi_M)$ is a complex manifold
$(M,J_M)$ with an exhausting plurisubharmonic function $\phi_M$.
Here plurisubharmonicity is always intended in the strict sense,
meaning that $-dd^c \phi_M = -d(d\phi_M \circ J_M)$ is a positive
$(1,1)$-form. We say that the Stein manifold is {\em complete} if
the gradient flow of $\phi_M$ exists for all positive times, and
{\em of finite type} if there is a $c_0$ such that all $c \geq c_0$
are regular values of $\phi_M$. Taking $\o_M = -dd^c\phi_M$ and
$\theta_M = -d^c\phi_M$ then makes $M$ into a convex symplectic
manifold, whose Liouville vector field is $\lambda_M = \nabla
\phi_M$ (and which therefore satisfies $d\phi_M(\lambda_M) > 0$ on
each regular level set of $\phi_M$). Completeness or finite type
nature of the Stein manifold imply the corresponding convex
symplectic properties.

\begin{lemma} \label{th:stein-rescale}
Let $(M,J_M,\phi_M)$ be a Stein manifold, and $h: \R \rightarrow \R$
a function satisfying $h'(c) > 0$, $h''(c) \geq 0$ for all $c$, and
such that there are $c_0$ and $\delta>0$ with $h''(c) \geq \delta
h'(c)$ for all $c \geq c_0$. Then $\tilde\phi_M = h(\phi_M)$ is
again plurisubharmonic and makes $(M,J_M)$ into a complete Stein
manifold. Denote the convex symplectic structure obtained from
$\tilde\phi_M$ by $(\tilde\o_M,\tilde\theta_M)$. If the original
Stein structure was complete and of finite type, there is a
diffeomorphism $f: M \rightarrow M$ such that $f^*(\tilde\theta_M) -
\theta_M =dR$ for some compactly supported function $R$.
\end{lemma}

\begin{lemma} \label{th:two-functions}
Let $\phi_M,\tilde\phi_M$ be two exhausting plurisubharmonic
functions on the same complex manifold $(M,J_M)$, of which the
second one is complete and of finite type (while the first one can
be arbitrary). Then there is an embedding $i: M \rightarrow M$ such
that $i^*(\tilde\theta_M) - \theta_M$ is an exact one-form.
\end{lemma}

Let $(J_{M,t},\phi_{M,t})$, $0 \leq t \leq 1$, be a smooth family of
Stein structures on a manifold $M$. We call this a {\em Stein
deformation} if the following two additional conditions hold: the
function $(t,x) \mapsto \phi_{M,t}(x)$ on $[0;1] \times M$ is
proper; and for each $t \in [0;1]$ there is a neighbourhood $t \in I
\subset [0;1]$ and a sequence $c_1 <c_2<\dots$ converging to
$+\infty$, such that $c_k$ is a regular level set for each
$\phi_{M,s}$, $s \in I$. A Stein deformation is called {\em
complete} if all the Stein structures in it are complete, and of
{\em finite type} if there is a $c_0$ such that all $c \geq c_0$ are
regular values of $\phi_{M,t}$ for all $t \in [0;1]$. Clearly, these
kinds of deformations induce the corresponding convex symplectic
notions.

It remains to make the connection with algebraic geometry. Let $X$
be a smooth projective variety, $E \rightarrow X$ an ample line
bundle, $s_E \in H^0(E)$ a nonzero holomorphic section, and $S =
s_E^{-1}(0)$ the hypersurface along which it vanishes. Ampleness
means that we can put a metric $||\cdot||_E$ on $E$ such that the
curvature form $\o_X = iF_{\nabla_E}$ of the associated connection
$\nabla_E$ is a positive $(1,1)$-form. The restriction of this form
to $M = X \setminus S$ can be written as $\o_M = \o_X|M =
-dd^c\phi_M$, where $\phi_M = -\log ||s_E||_E$. This is clearly an
exhausting function, hence defines a Stein structure.

\begin{lemma} \label{th:normal-crossing}
Suppose that $S$ has only normal crossing singularities (but $s_E$
can vanish along the irreducible components with arbitrary
multiplicities). Then $\phi_M$ is of finite type.
\end{lemma}

There is also a version of this for deformations: in the same
algebro-geometric situation, given a family $||\cdot||_{E,t}$ of
metrics, one gets a finite type Stein deformation $(J_{M,t} = J_M,
\phi_{M,t} = -\log ||s_E||_{E,t})$.

\section{Stein-Essential Lagrangian submanifolds\label{sec:essential}}

Following a line of thought similar to the one in
\cite{biran-cieliebak01}, we combine ``soft'' displacement methods
for subcritical Stein manifolds with ``hard'' Lagrangian
intersection results to derive some restrictions on embeddings of
Stein manifolds.

Let $(M,\o_M,\theta_M)$ be any exact symplectic manifold. By a
Hamiltonian isotopy of $M$, we will mean an isotopy $(g_t)$, $0 \leq
t \leq 1$, starting with $g_0 = id$, which is induced by a smooth
family of Hamiltonian functions $H_t$, such that $(t,x) \mapsto
H_t(x)$ has compact support in $[0;1] \times M$. The {\em Hofer
length} of $(g_t)$ is defined as \[\int_0^1 \big( \max(H_t) -
\min(H_t) \big)\, dt.\] Now let $L \subset M$ be a Lagrangian
submanifold (throughout, all such submanifolds will be assumed to be
compact). Consider Hamiltonian isotopies $(g_t)$ such that $g_t(L)
\cap L = \emptyset$. The infinum of the Hofer lengths of all these
isotopies is called the {\em displacement energy} of $L$ (of course,
there are cases where no such isotopy exists, and then the
displacement energy is $\infty$). Recall that by definition, a
Lagrangian isotopy $(L_t)$ is {\em exact} if the class
$[\theta_M|L_t] \in H^1(L_t;\R) \iso H^1(L_0;\R)$ is constant in
$t$. These are precisely the Lagrangian isotopies which can be
embedded into Hamiltonian ones, in the sense that there is a $(g_t)$
with $g_t(L_0) = L_t$. Hence, the displacement energy is invariant
under exact Lagrangian isotopies. Chekanov's theorem
\cite{chekanov95,chekanov98,oh97b} says:

\begin{theorem} \label{th:chekanov}
Let $L \subset M$ be a compact Lagrangian submanifold whose
displacement energy is $E < \infty$. Let $J_M$ be an
$\o_M$-compatible almost complex structure which is convex at
infinity. Then there is a non-constant $J_M$-holo\-mor\-phic map $u:
(\D,\partial \D) \rightarrow (M,L)$, where $\D \subset \C$ is the
closed unit disc, whose area is $\int u^*\o_M \leq E$. \qed
\end{theorem}

Convexity at infinity of the almost complex structure means that
there is an exhausting function $\phi_M$ such that outside a compact
subset, $-d(d\phi_M \circ J_M)$ is positive on all $J_M$-complex
tangent planes. This holds for the given complex structure on any
Stein manifold, but it also allows one to deform that structure
(compatibly with the symplectic form) on a compact subset.
Chekanov's theorem actually holds in somewhat greater generality,
but that will not be necessary for our purpose.

Recall that a Stein manifold $(M,J_M,\phi_M)$ is {\em subcritical}
if $\phi_M$ is a Morse function and has only critical points of
index $< \frac{1}{2}\dim_\R M$. The next statement is a special case
of \cite[Lemma 3.2]{biran-cieliebak01} (and technically somewhat
simpler than the general result):

\begin{lemma} \label{th:subcritical}
For any Lagrangian submanifold $L$ in a complete subcritical Stein
manifold $M$, there is an exact Lagrangian isotopy $(L_t)$ such that
$L_0 = L$, $L_1 \cap L_0 = \emptyset$. \qed
\end{lemma}

Suppose that we have a Stein manifold, containing a Lagrangian
submanifold such that there are {\em no} non-constant holomorphic
discs bounding it. Then our Stein manifold cannot be subcritical
(one would use Lemma \ref{th:stein-rescale} to make it complete, and
then combine Theorem \ref{th:chekanov} with Lemma
\ref{th:subcritical}). A little less obviously, this manifold cannot
have an exact symplectic embedding into any subcritical Stein
manifold. We will spend the rest of this section deriving a more
complicated version of the latter statement.
%
% note: 1. make M complete, leaving L Lagrangian
% 2. choose c_0 > 0 such that L is contained in the corresponding
% sublevel set. We can replace M by the open subset obtained by
% taking the cone map (0;\infty \times phi^{-1}(c_0) \rightarrow
% \phi^{-1}([c_0;\infty)). This is still Stein, it just excludes
% other handles which are attached - and it still has the same
% property, of not having holomorphic discs - the upshot is,
% we can also assume M to be of finite type
% 3. now go as in the lemma afterwards, taking the deformation
% to be constant...

Let $(M,J_M,\phi_M)$ be a complete finite type Stein manifold, and
take $c_0$ as usual. On $\phi_M^{-1}([c_0;\infty))$, which in
symplectic terms is the cone part of $M$, consider the splitting
\begin{equation} \label{eq:perp}
TM = \xi_M \oplus \xi_M^\perp
\end{equation}
where $\xi_M = \ker(d\phi_M) \cap \ker(d^c\phi_M)$ is the contact
hyperplane field on each level set $\phi_M^{-1}(c)$, $c \geq c_0$;
and $\xi_M^\perp = \R\nabla\phi_M \oplus \R\rho_M$ is spanned by the
Liouville vector field together with the Reeb vector field on each
level set, which is $\rho_M = J_M\nabla\phi_M/||\nabla \phi_M||^2$.
The decomposition \eqref{eq:perp} is $J_M$-invariant and orthogonal
with respect to $\o_M$. One can therefore find an $\o_M$-compatible
almost complex structure $\tilde{J}_M$ which
\begin{equation} \label{eq:tildej}
\begin{cases}
\text{is equal to $J_M$}
& \text{on $\phi_M^{-1}(-\infty;c_0])$}; \\
\text{preserves $\xi_M$, and maps $\lambda_M$ to} & \\
\qquad \text{a positive multiple of $\rho_M$}
& \text{on $\phi_M^{-1}([c_0;+\infty))$}; \\
\text{is invariant under the Liouville flow} & \text{on
$\phi_M^{-1}([c_0+1;+\infty))$.}
\end{cases}
\end{equation}

\begin{lemma} \label{th:maximum}
Let $\Sigma_0$ be a compact connected Riemann surface with boundary,
and $u: \Sigma_0 \rightarrow M$ a $\tilde{J}_M$-holomorphic map such
that $u(\partial\Sigma_0) \subset \phi_M^{-1}((-\infty;c_0])$. Then
$u(\Sigma_0) \subset \phi_M^{-1}((-\infty;c_0])$.
\end{lemma}

\proof The second part of \eqref{eq:tildej} implies that on
$\phi_M^{-1}([c_0;+\infty))$, $-d\phi_M \circ \tilde{J}_M = \eta
\,\theta_M$ with a strictly positive function $\eta$. Hence
\begin{equation} \label{eq:almost-laplace}
 -d(d\phi_M \circ \tilde{J}_M) = \eta \, \o_M +
 \frac{d\eta}{\eta} \wedge (-d\phi_M \circ \tilde{J}_M).
\end{equation}
Suppose that we have a $\tilde{J}_M$-holomorphic map $u: \Sigma_0
\rightarrow M$ with $u(\partial\Sigma_0) \subset \phi_M^{-1}
((-\infty;c_0])$ but $u(\Sigma_0) \not\subset \phi_M^{-1}
((-\infty;c_0])$. Then there is a $c > c_0$ such that $u$ intersects
the level set $\phi_M^{-1}(c)$ transversally in a nonempty set.
Consider the function $\psi = \phi_M \circ u$ on the surface $\Sigma
= u^{-1}\phi_M^{-1}([c;+\infty)) \subset \Sigma_0$. By pulling back
\eqref{eq:almost-laplace} and using the positivity of $u^*\o_M$ we
obtain a differential inequality for $\psi$, which in a local
holomorphic coordinate $z = s+it$ can be written as
\[
 (\partial_s^2 + \partial_t^2) \psi
 - \sigma(s,t) \partial_s \psi - \tau(s,t) \partial_t \psi \geq 0
 \]
with $\sigma = (\eta \circ u)^{-1} \partial_s (\eta \circ u)$, $\tau
= (\eta \circ u)^{-1} \partial_t (\eta \circ u)$. The strong maximum
principle \cite[Theorem 3.5]{gilbarg-trudinger} applies to solutions
of such equations, hence $\psi \leq c$ everywhere on $\Sigma$, which
means that $u(\Sigma_0) \subset \phi_M^{-1}((-\infty;c])$. Since $c$
can be chosen arbitrarily close to $c_0$, the result follows. \qed

Because $\tilde{J}_M$ is invariant under the Liouville flow outside
a compact subset, it is tame in the sense of \cite[Chapter V,
Definition 4.1.1]{audin-lafontaine}. In particular, the monotonicity
lemma \cite[Chapter V, Proposition 4.3.1]{audin-lafontaine} applies:

\begin{lemma} \label{th:mono}
Let $\tilde g_M$ be the metric associated to $\o_M$ and
$\tilde{J}_M$. There is a $\rho>0$, which is less than the
injectivity radius of $\tilde g_M$, and an $\epsilon>0$, such that
the following holds. Let $x$ be any point in $M$, and $B_r(x)$ the
closed ball of radius $r \leq \rho$ around it. If $\Sigma$ is a
compact Riemann surface with boundary, and $u: \Sigma \rightarrow
B_r(x)$ a $\tilde{J}_M$-holomorphic map satisfying $x \in u(\Sigma)$
and $u(\partial\Sigma) \subset
\partial B_r(x)$, then $\int u^*\o_M \geq \epsilon\, r^2$. \qed
\end{lemma}

\begin{lemma} \label{th:energy}
For every $E>0$ there is a $C>0$ with the following property. Let
$\Sigma$ be a compact connected Riemann surface with boundary, whose
boundary is decomposed into two nonempty unions of circles
$\partial_-\Sigma \cup \partial_+\Sigma$. Let $u: \Sigma \rightarrow
M$ be a $\tilde{J}_M$-holomorphic map such that $u(\partial_-\Sigma)
\subset \phi_M^{-1}((-\infty;c_0])$ and $u(\partial_+\Sigma) \subset
\phi_M^{-1}([C;\infty))$. Then $\int_{\Sigma} u^*\o_M > E$.
\end{lemma}

\proof Consider the diffeomorphism $f: [0;\infty) \times
\phi_M^{-1}(c_0) \rightarrow \phi_M^{-1}([c_0;\infty))$ which
defines the conical end structure. Since the metric $\tilde g_M$
blows up on the cone, the distance between any two sets $f(\{i\}
\times \phi_M^{-1}(c_0))$, $i = 0,1,2,\cdots$ is bounded below by
some $\delta>0$. Take the constants $\rho,\epsilon$ from Lemma
\ref{th:mono}. After making $\delta$ smaller if necessary, we may
assume that $\delta/2<\rho$; we then take an integer $k$ greater
than $9 \delta^{-2} \epsilon^{-1} E$, and choose $C$ so that
$\phi_M^{-1}([C;\infty)) \subset f([k;\infty) \times
\phi_M^{-1}(c_0))$.

Since $\Sigma$ is connected and intersects both
$\phi_M^{-1}((-\infty;c_0])$ and $f([k;\infty) \times
\phi_M^{-1}(c_0))$ nontrivially, there are points $z_1,\dots,z_k \in
\Sigma$ such that $x_i = u(z_i) \in f(\{i-1/2\} \times
\phi_M^{-1}(c_0))$. The balls $B_r(x_i)$, for any $r < \delta/2$,
are mutually disjoint. Choose $r \in (\delta/3;\delta/2)$ in such a
way that $u$ is transverse to all the boundaries $\partial
B_r(x_i)$. By Lemma \ref{th:mono}, each $u_i = u|u^{-1}(B_r(x_i)):
u^{-1}(B_r(x_i)) \rightarrow B_r(x_i)$ has area $\geq \epsilon
\delta^2/9$. Hence, the total area of $u$ is $\geq k \epsilon
\delta^2/9 > E$. \qed

Let $(M,J_M,\phi_M)$ be a finite type Stein manifold. We say that a
compact Lagrangian submanifold $L \subset M$ is {\em
Stein-essential} if for each $E>0$ there is a finite type Stein
deformation $(J_{M,t}, \phi_{M,t})$ and a smooth family of compact
submanifolds $L_t \subset M$ ($0 \leq t \leq 1$), with the following
properties: at the starting point $t = 0$, we have the original
Stein structure and Lagrangian submanifold $L = L_0$; for all $t$,
$L_t$ is $\o_{M,t}$-Lagrangian, and the cohomology class
$[\theta_{M,t}] \in H^1(L_t;\R) \iso H^1(L;\R)$ is constant in $t$;
and at the opposite end, every $J_{M,1}$-holomorphic map $u:
(\D,\partial \D) \rightarrow (M,L_1)$ with $\int u^*\o_{M,1} \leq E$
is constant.

\begin{prop} \label{th:proposition}
Let $M$ be a finite type Stein manifold, which admits an embedding
$i: M \rightarrow N$ into a complete subcritical Stein manifold,
such that $i^*\theta_N - \theta_M$ is an exact one-form. Then $M$
cannot contain any Stein-essential Lagrangian submanifolds.
\end{prop}

\proof Assume that on the contrary, there is an Stein-essential
Lagrangian submanifold $L \subset M$. By definition, for each $E$ we
can find a finite type Stein deformation $(J_{M,t},\phi_{M,t})$ and
family of Lagrangian submanifolds $(L_t)$, such that $E$ is a strict
lower bound for the area of non-constant $J_{M,1}$-holomorphic discs
in $(M,L_1)$. We may take $E$ to be the displacement energy of
$i(L)$ inside $N$, which is finite by Lemma \ref{th:subcritical}.

Because the deformation $(J_{M,t},\phi_{M,t})$ is of finite type,
there is a $c_0>0$ so that all $c \geq c_0$ are regular values of
$\phi_{M,t}$ for all $t$. After making $c_0$ larger if necessary,
one can also assume that $L_t \subset \phi_{M,t}^{-1}((-\infty;
c_0])$. Lemma \ref{th:stein-rescale} says that one can find
functions $h_t$ depending smoothly on $t$, such that the modified
Stein structures $(J_{M,t}, \tilde\phi_{M,t} = h_t(\phi_{M,t}))$ are
complete. Choose these functions in such a way that $h_t(c) = c$ for
$c \leq c_0$, which means that the $L_t$ remain Lagrangian for the
associated modified convex symplectic structures
$(\tilde\o_{M,t},\tilde\theta_{M,t})$.

By construction $(J_{M,1},\tilde\phi_{M,1})$ is complete and of
finite type. Introduce a new $\tilde\o_{M,1}$-compatible almost
complex structure $\tilde{J}_{M,1}$ as in the discussion preceding
Lemma \ref{th:maximum}. More explicitly, to carry over that
construction to the current situation, one should replace the
notation $J_M$, $\phi_M$, $\tilde{J}_M$ in \eqref{eq:tildej} by
$J_{M,1}$, $\tilde\phi_{M,1}$, $\tilde{J}_{M,1}$ respectively; and
similarly $\xi_M$,$\lambda_M$,$\rho_M$ are now the contact
hyperplane field, Liouville vector field, and Reeb vector field
associated to $(\tilde\o_{M,1}, \tilde\theta_{M,1})$ and to the
conical end $[c_0;+\infty) \times \tilde\phi_{M,1}^{-1}(c_0)
\rightarrow \tilde\phi_{M,1}^{-1}([c_0;+\infty))$. We will now state
some properties of the data introduced so far:

\begin{itemize}
\item[\bf (a)] {\em There is an embedding $\tildei: M \rightarrow
N$ with $\tildei(L) = i(L)$, such that $\tildei^*\theta_N -
\tilde{\theta}_{M,0}$ is an exact one-form.}
\end{itemize}

To obtain that, restrict $i$ to an embedding of
$\phi_M^{-1}((-\infty;c_0])$ into $N$, note that $\theta_M =
\tilde{\theta}_{M,0}$ on that subset, and then extend it to the
whole of $(M,\tilde{\o}_{M,0}, \tilde\theta_{M,0})$ using Lemma
\ref{th:weinstein-embed}.

\begin{itemize}
\item[{\bf (b)}] {\em There is a diffeomorphism $\tilde{f}_1: M
\rightarrow M$ such that $\tilde{f}_1^*\tilde\theta_{M,1} -
\tilde\theta_{M,0}$ is an exact one-form, and $\tilde{f}_1(L_0)$ is
exact Lagrangian isotopic to $L_1$.}
\end{itemize}

By definition $(J_{M,t},\tilde{\phi}_{M,t})$ is a complete finite
type deformation, so Lemma \ref{th:weinstein-deform-2} provides a
family of diffeomorphisms $\tilde{f}_t: M \rightarrow M$ such that
$\tilde{f}_t^*\tilde\theta_{M,t} - \tilde\theta_{M,0}$ are exact.
$\tilde{f}_1\tilde{f}_t^{-1}(L_t)$ is a Lagrangian isotopy between
$\tilde{f}_1(L_0)$ and $L_1$, and the cohomology class
$[\tilde{\theta}_{M,1}|\tilde{f}_1\tilde{f_t}^{-1}(L_t)] =
[\tilde{\theta}_{M,0}|\tilde{f_t}^{-1}(L_t)] =
[\tilde\theta_{M,t}|L_t] = [\theta_{M,t}|L_t]$ is constant in $t$,
which means that the isotopy is exact.

\begin{itemize}
\item[{\bf (c)}] {\em The image of any
$\tilde{J}_{M,1}$-holomorphic disc $\tilde{u}: (\D,\partial \D)
\rightarrow (M,L_1)$ lies in
$\tilde\phi_{M,1}^{-1}((-\infty;c_0])$}.
\end{itemize}

By construction $L_1 \subset \tilde\phi_{M,1}^{-1}((-\infty;c_0])$;
therefore Lemma \ref{th:maximum} applies and yields the desired
result.

\begin{itemize}
\item[{\bf (d)}] {\em There is a $C>c_0$ with the following
property. Let $\Sigma$ be a compact connected Riemann surface with
boundary, whose boundary is decomposed into two nonempty unions of
circles $\partial_-\Sigma \cup \partial_+\Sigma$. Let $\tilde{u}:
\Sigma \rightarrow M$ be a $\tilde{J}_{M,1}$-holomorphic map such
that $\tilde{u}(\partial_-\Sigma) \subset
\tilde\phi_{M,1}^{-1}((-\infty;c_0])$ and
$\tilde{u}(\partial_+\Sigma) \subset
\tilde\phi_{M,1}^{-1}([C;\infty))$. Then $\int_{\Sigma}
u^*\tilde{\o}_{M,1} > E$.}
\end{itemize}

Up to the change in notation, this is Lemma \ref{th:energy}.

Consider the compact subset $K =
\tilde\phi_{M,1}^{-1}((-\infty;C+1])$, and let $j = \tildei \circ
\tilde{f}_1^{-1}|K : K \rightarrow N$. Clearly, one can find an
$\o_N$-compatible almost complex structure $\tilde{J}_N$ with the
following two properties: $\tilde{J}_N = J_N$ outside a compact
subset; and $j^*\tilde{J}_N = \tilde{J}_{M,1}|K$. From (b) we know
that the Lagrangian submanifold $j(L_1)$ is exact isotopic to
$\tildei(L)$, hence its displacement energy is again $E$. Since
$\tilde{J}_N = J_N$ at infinity, we can apply Theorem
\ref{th:chekanov}, which shows that there is a non-constant
$\tilde{J}_N$-holomorphic disc $u: (\D,\partial \D) \rightarrow
(N,j(L_1))$ with $\int u^*\o_N \leq E$.

Choose a $c \in [C;C+1]$ such that $u$ intersects the hypersurface
$j(\tilde\phi_{M,1}^{-1}(c))$ transversally. Consider only the part
of our $\tilde{J}_N$-holomorphic disc $u$ that lies on the interior
side of that hypersurface. This may have several connected
components; we ignore all of them except the one which contains
$\partial \D$, and compose that with $j^{-1}$ to obtain a
$\tilde{J}_{M,1}$-holomorphic map $\tilde{u}: \Sigma \rightarrow K
\subset M$. By construction, $\Sigma$ is a connected compact Riemann
surface with boundary; its boundary contains one circle
$\partial_-\Sigma$ such that $\tilde{u}(\partial_-\Sigma) \subset
L_1 \subset \tilde\phi_{M,1}^{-1}((-\infty;c_0])$, and if
$\partial_+\Sigma$ is the union of all the other boundary circles,
then $\tilde{u}(\partial_+\Sigma) \subset \tilde\phi_{M,1}^{-1}(c)
\subset \tilde\phi_{M,1}^{-1}([C;\infty))$; finally $\int
\tilde{u}^*\tilde\o_{M,1} \leq E$. By (d) above, this is possible
only if $\partial_+\Sigma = \emptyset$, which means that $\tilde{u}$
is a non-constant $\tilde{J}_{M,1}$-holomorphic disc in $(M,L)$.
Applying (c) we find that the image of $\tilde{u}$ must be contained
in $\tilde\phi_{M,1}^{-1}((-\infty;c_0])$, which implies that it is
in fact a $J_{M,1}$-holomorphic disc, with $\int \tilde{u}^*\o_{M,1}
\leq E$. However, given our original choice of the deformation, the
existence of such a disc violates the definition of Stein-essential
Lagrangian submanifold. \qed

In fact, the requirement that $N$ is complete can be omitted, due to
the following observation, which is similar to step (a) in the
previous proof:

\begin{lemma} \label{th:stupid}
Let $M$ be a finite type Stein manifold. If $M$ admits an embedding
$i: M \rightarrow N$ into a subcritical Stein manifold, such that
$i^*\theta_N - \theta_M$ is an exact one-form, then it also admits
an embedding into a complete subcritical Stein manifold, with the
same property.
\end{lemma}

\proof Take $c_0$ so that all $c \geq c_0$ are regular values of
$\phi_M$. Use Lemma \ref{th:stein-rescale} to find an $h$ such that
$\tilde\phi_N = h(\phi_N)$ gives rise to a complete Stein structure.
This is still subcritical, because the critical points and their
Morse indices remain the same. $h$ can be chosen in such a way that
the new convex symplectic structure $(\tilde\o_N,\tilde\theta_N)$
agrees with the old one on $i(\phi_M^{-1}((-\infty;c_0]) \subset N$.
By restricting $i$ to $\phi_M^{-1}((-\infty;c_0])$, and then
extending it again using Lemma \ref{th:weinstein-embed}, one gets an
embedding $j: M \rightarrow N$ such that $j^*\tilde\theta_N -
\theta_M$ is exact. \qed

\section{Linking tori}\label{sec:linking}

Throughout this section, $X$ will be a smooth projective algebraic
surface; $S \subset X$ an algebraic curve with only normal crossing
singularities; and $p \in S$ a crossing point. Set $M = X \setminus
S$. Take local holomorphic coordinates $(a,b)$ centered at $p$ in
which $S = \{ab = 0\}$, and let $U \subset X$ be a ball of some
radius $\rho>0$ in those coordinates. Consider the torus $L \subset
M$ given by $\{|a| = \mu, \; |b| = \nu\}$ for some $0 <\mu,\nu <
\rho/\sqrt{2}$. We will call such an $L$, as well as any other torus
isotopic to it inside $U \cap M$, a {\em linking torus} for $S$ at
$p$.

Recall that, given any algebraic curve on a smooth algebraic
surface, one can resolve its singularities by blowups, until only
normal crossings remain (\cite{BPV}, Chapter II). The linking tori
constructed in this way can be viewed as lying in the complement of
the original curve, since blowups leave that complement unchanged.
We will now consider in more detail the simplest example of this,
which is relevant for our application later on. Let $\bar{X}$ be a
smooth projective algebraic surface, and $\bar{S} \subset \bar{X}$ a
curve which has a cusp singularity at the point $\bar{p}$. Blow up
to get a map $q: X \rightarrow \bar{X}$, such that $S =
q^{-1}(\bar{S})$ has only normal crossings. We assume that this
resolution is the minimal one (meaning that no exceptional component
of $S$ can be blown down without violating the normal crossing
condition). Take local coordinates $(c,d)$ centered at $\bar{p}$ in
which $\bar{S} = \{c^2 = d^3\}$; let $\bar{V} \subset \bar{X}$ be a
small ball in these coordinates; and set $V = q^{-1}(\bar{V})$. The
curve $S \cap V$ consists of a small piece of the principal
component, which is the proper transform of $\bar{S}$, and three
exceptional components of multiplicities $2$, $3$ and $6$. Figure
\ref{fig:cusp} summarizes the stages of the blowup process and the
corresponding coordinate changes (the thick lines are the
exceptional components, and the dots indicate the origin of the
coordinate systems used).

\begin{figure}
\begin{center}
\begin{picture}(0,0)%
\includegraphics{cusp.pstex}%
\end{picture}%
\setlength{\unitlength}{3947sp}%
\begingroup\makeatletter\ifx\SetFigFont\undefined%
\gdef\SetFigFont#1#2#3#4#5{%
  \reset@font\fontsize{#1}{#2pt}%
  \fontfamily{#3}\fontseries{#4}\fontshape{#5}%
  \selectfont}%
\fi\endgroup%
\begin{picture}(5196,4684)(504,-4358)
\put(701,-3770){\makebox(0,0)[lb]{\smash{{\SetFigFont{12}{14.4}{\rmdefault}{\mddefault}{\updefault}{\color[rgb]{0,0,0}2}%
}}}}
\put(1801,-136){\makebox(0,0)[lb]{\smash{{\SetFigFont{12}{14.4}{\rmdefault}{\mddefault}{\updefault}{\color[rgb]{0,0,0}$c^2 - d^3$ }%
}}}}
\put(1801,-1261){\makebox(0,0)[lb]{\smash{{\SetFigFont{12}{14.4}{\rmdefault}{\mddefault}{\updefault}{\color[rgb]{0,0,0}$t^2(s^2-t)$}%
}}}}
\put(1801,-2386){\makebox(0,0)[lb]{\smash{{\SetFigFont{12}{14.4}{\rmdefault}{\mddefault}{\updefault}{\color[rgb]{0,0,0}$u^2v^3(v-u)$}%
}}}}
\put(2551,-1711){\makebox(0,0)[lb]{\smash{{\SetFigFont{12}{14.4}{\rmdefault}{\mddefault}{\updefault}{\color[rgb]{0,0,0}$t = uv$}%
}}}}
\put(2551,-1936){\makebox(0,0)[lb]{\smash{{\SetFigFont{12}{14.4}{\rmdefault}{\mddefault}{\updefault}{\color[rgb]{0,0,0}$s = v$}%
}}}}
\put(2551,-511){\makebox(0,0)[lb]{\smash{{\SetFigFont{12}{14.4}{\rmdefault}{\mddefault}{\updefault}{\color[rgb]{0,0,0}$c = st$}%
}}}}
\put(2551,-736){\makebox(0,0)[lb]{\smash{{\SetFigFont{12}{14.4}{\rmdefault}{\mddefault}{\updefault}{\color[rgb]{0,0,0}$d = t$}%
}}}}
\put(1801,-3661){\makebox(0,0)[lb]{\smash{{\SetFigFont{12}{14.4}{\rmdefault}{\mddefault}{\updefault}{\color[rgb]{0,0,0}$-a^6b(b+1)^2$}%
}}}}
\put(2551,-2911){\makebox(0,0)[lb]{\smash{{\SetFigFont{12}{14.4}{\rmdefault}{\mddefault}{\updefault}{\color[rgb]{0,0,0}$v = a$}%
}}}}
\put(2551,-3136){\makebox(0,0)[lb]{\smash{{\SetFigFont{12}{14.4}{\rmdefault}{\mddefault}{\updefault}{\color[rgb]{0,0,0}$u = (b+1)a$}%
}}}}
\put(4501,-1786){\makebox(0,0)[lb]{\smash{{\SetFigFont{12}{14.4}{\rmdefault}{\mddefault}{\updefault}{\color[rgb]{0,0,0}$d = (b+1)a^2$}%
}}}}
\put(4501,-1561){\makebox(0,0)[lb]{\smash{{\SetFigFont{12}{14.4}{\rmdefault}{\mddefault}{\updefault}{\color[rgb]{0,0,0}$c = (b+1)a^3$}%
}}}}
\put(701,-3211){\makebox(0,0)[lb]{\smash{{\SetFigFont{12}{14.4}{\rmdefault}{\mddefault}{\updefault}{\color[rgb]{0,0,0}3}%
}}}}
\end{picture}%
\caption{\label{fig:cusp}}
\end{center}
\end{figure}

Consider the linking torus $L = \{|a| = \mu, |b| = \nu\}$ at the
point $(a,b) = (0,0)$ where the principal component of $S \cap V$
crosses the exceptional component of multiplicity $6$ (this point is
indicated by the small arrow in Figure \ref{fig:cusp}). Its image
under $q$ is the torus $\bar{L} = \bar{L}_{\mu,\nu}$ parametrized by
\[
 c = \mu^3(\nu e^{i\gamma} + 1)e^{3i\delta}, \quad
 d = \mu^2(\nu e^{i\gamma} + 1)e^{2i\delta}
\]
for $(\gamma,\delta) \in \R/2\pi\Z$, and where $\mu,\nu > 0$ are
suitably small constants. By keeping $\mu$ constant and letting $\nu
\rightarrow 0$, one obtains a smooth family of tori in $\bar{V}
\setminus \bar{S}$, which in the limit shrink to the loop $(c =
\mu^3 e^{3i\delta}, d = \mu^2 e^{2i\delta})$ lying on $(\bar{S}
\setminus \{(0,0)\}) \cap \bar{V}$.

The topological aspect of cusp singularities is well-known: the
intersection $S^3 \cap \{c^2=d^3\}$ is a $(2,3)$-torus knot, which
is a trefoil $\kappa$. One can find a diffeomorphism $\bar{V}
\setminus \{(0,0)\} \iso (0;1) \times S^3$ which takes $\bar{S}
\setminus \{(0,0)\}$ to $(0;1) \times \kappa$, hence identifies
$\bar{V} \setminus \bar{S}$ with $(0;1) \times (S^3 \setminus
\kappa)$. From the argument given above, it follows that the loop on
$\bar{L}$ given by $\{\gamma = \text{\it const.}\}$ is homotopic to
a longitude of $\kappa$. Here, by longitude we mean a curve in $S^3
\setminus \kappa$ which runs parallel to $\kappa$, for some framing
which may not necessarily be the canonical one (this ambiguity could
be settled by explicit computation, but it is irrelevant for our
purpose). Similarly, by inspection of the limit $\nu \rightarrow 0$
with fixed $\mu$ and $\delta$, one sees that the other loop
$\{\delta = \text{\it const.}\}$ on $\bar{L}$ is a meridian of
$\kappa$. It is a general fact about nontrivial knots (Dehn's Lemma,
see e.g.\ \cite[Theorem 11.2]{lickorish}) that longitude and
meridian together define an injective homomorphism $\Z^2 \rightarrow
\pi_1(S^3 \setminus \kappa)$, which for us means that
$\pi_1(\bar{L}) \rightarrow \pi_1(\bar{V} \setminus \bar{S})$ is
injective. Using the identification $V \setminus S \iso \bar{V}
\setminus \bar{S}$ provided by $q$, we arrive at this conclusion:

\begin{lemma} \label{th:trefoil}
If a loop on $L$ bounds a disc in $V \setminus S$, then it must be
contractible on $L$ itself. \qed
\end{lemma}

Returning to the general discussion of linking tori, we now
reformulate their definition using a degeneration of $X$ to a normal
crossing surface. Let $Y$ be the variety obtained by blowing up
$(0,p) \in \P^1 \times X$, and $\pi: Y \rightarrow \P^1 \times X
\rightarrow \P^1$ the projection to the first variable. The smooth
fibres $Y_t$, $t \in \P^1 \setminus \{0\}$, are obviously isomorphic
to $X$. The singular fibre has two irreducible components, $Y_0 = Z
\cup P$: $Z$ is the blowup of $X$ at $p$, and $P = {\mathbb P}(\C
\oplus TX_p)$ is the exceptional divisor in $Y$. They are joined
together by a normal crossing, where one identifies the exceptional
curve in $Z$ with the line $C_0 = \P(\{0\} \oplus TX_p)$ in $P$. Let
$T \subset Y$ be the proper transform of $\P^1 \times S \subset \P^1
\times X$ under the blowup, and $T_t = T \cap Y_t$. For $t \neq 0$
one can obviously identify $T_t \subset Y_t$ with $S \subset X$;
while $T_0$ is the union of $T \cap Z$, which is the proper
transform of $S$ under blowing up $p \in X$, and of $T \cap P$,
which consists of the two lines $C_1,C_2 \subset P$ obtained by
projectivizing $\C \times \text{\it (tangent space to either branch
of $S$ at $p$)}$. Choose an isomorphism $P \iso \P^2$ in such a way
that the $C_k$ become the coordinate lines, and let $K_0 \subset
(\C^*)^2 \iso P \setminus (C_0 \cup C_1 \cup C_2)$ be one of the
standard Clifford tori. One can find a submanifold with boundary $K
\subset Y$, lying in $\pi^{-1}([-1;1])$ and with boundary in
$\pi^{-1}(\{-1;1\})$, such that $\pi|K: K \rightarrow [-1;1]$ is a
smooth fibration, whose fibre over $t = 0$ is the given $K_0$. One
way to think of this is to choose a connection (a horizontal
subbundle) on the open subset of $\pi$-regular points of $Y$. Since
$K_0$ lies in that subset, one can use parallel transport to move it
to other fibres, and doing that in both directions along the real
axis yields $K$. In fact, any $K$ with the properties stated above
can be obtained in this way, for some choice of connection.

\begin{lemma} \label{th:linking-torus}
For all sufficiently small $t \in [-1;1] \setminus \{0\}$, $K_t = K
\cap Y_t$ is a linking torus for $S$ at the crossing point $p$.
\end{lemma}

It may be appropriate to first clarify the meaning of this. As
before, let $(a,b)$ be coordinates centered at $p$ in which $S =
\{ab = 0\}$, and $U$ a ball in those coordinates. By identifying
$Y_t \iso X$ for $t \in [-1;1] \setminus \{0\}$, one can think of
the $K_t$ as tori inside $X$. A more technical formulation of the
Lemma is that for sufficiently small such $t$, $K_t$ lies in $U
\setminus S$, and is isotopic inside that set to the standard
linking torus $\{|a| = \mu, \; |b| = \nu\}$.  (This formulation, and
the following proof, are somewhat pedantic, but are engineered to
adapt well to the symplectic geometry requirements to be imposed
subsequently.)

\proof Let $W$ be the preimage of $\P^1 \times U$ under the blowup
map $Y \rightarrow \P^1 \times X$. Since $K_0 \subset P \subset W$,
one has $K_t \subset W \cap Y_t = U$ for sufficiently small $t \neq
0$. Similarly, because $K_0 \cap T = K_0 \cap (C_1 \cup C_2) =
\emptyset$, one has $K_t \cap S = \emptyset$ for sufficiently small
$t \neq 0$. The next step is to show that the isotopy type of $K_t$
inside $U \setminus S$ is independent of the choice of $K$. If one
thinks of that choice as given by a connection, any two connections
can be deformed into each other, which gives rise to an isotopy of
the associated submanifolds $K_t$. The previous considerations show
that for small $t$, this isotopy will take place inside $U \setminus
S$.

With that in mind, it is sufficient to prove the statement that the
$K_t$ are linking tori for just one choice of $K$. We write down the
local picture near $P \subset Y$ in coordinates:
\[
\begin{aligned}
 & Y = \{ (t,a,b,[\tau:\alpha:\beta]) \in \C^3 \times \P^2 \suchthat
 (t,a,b) \in [\tau:\alpha:\beta]\}, \\
 & \pi(t,a,b,[\tau:\alpha:\beta]) = t, \\
 & Z = \{ t = 0,\; \tau = 0 \}, \\
 & P = \{ t = a = b = 0 \}, \\
 & T = \{ \alpha\beta = 0 \}.
\end{aligned}
\]
Here $(t,a,b)$ should actually lie in a small neighbourhood of
$(0,0,0)$, but we omit that to make the notation more transparent.
By definition, $K_0 = \{\tau = 1, \; |\alpha| = \mu, \; |\beta| =
\nu \}$ for some constants $\mu,\nu > 0$, and one can therefore take
$K = \{t \in \R, \; |a| = \mu |t|, \; |b| = \nu |t|\}$, in which
case $K_t$ ($t \neq 0$) is clearly a family of linking tori, whose
diameter shrinks as $t \rightarrow 0$. \qed

From this point onwards, we will make the additional assumption that
there is an effective divisor $D$ on $X$ whose support is $S$ (in
other words, $D$ is a sum of the irreducible components of $S$ with
positive multiplicities), which is ample. Let $E = \O_X(D)$ be the
associated ample line bundle. Form the tensor product $\O_{\P^1}(1)
\boxtimes E \rightarrow \P^1 \times X$ and pull it back to a line
bundle on $Y$ (keeping the notation for simplicity). For $d \gg 0$,
$F = (\O_{\P^1}(1) \boxtimes E)^{\otimes d} \otimes \O_Y(-P)
\rightarrow Y$ is again ample. We will now recall Kodaira's
classical proof of this fact; for full details see e.g.\ \cite[p.\
185]{griffiths-harris}. Start with a metric $||\cdot||_E$ on $E$
whose curvature (more precisely $iF_{\nabla_E}$, where $\nabla_E$ is
the associated connection) is a positive $(1,1)$-form, denoted by
$\o_X$. Similarly, on $\O_{\P^1}(1)$ we choose a metric whose
curvature is a positive $(1,1)$-form $\o_{\P^1}$. Tensor them
together to give a metric on $\O_{\P^1}(1) \boxtimes E$, with
curvature $\o_{\P^1} + \o_X$. By specializing to the point $(0,p)
\in \P^1 \times X$, this induces a metric on the bundle $\O_P(1)
\rightarrow P$, and a Fubini-Study form $\o_P$ on $P$. One can
identify $\O_Y(-P)|P \iso \O_P(1)$, so this gives a metric on
$\O_Y(-P)|P$, which one can extend to a small neighbourhood of $P$.
On the complement of $P$, $\O_Y(-P)$ is canonically trivial, so one
can take a constant metric, and patch that together with the other
one using a cutoff function. The outcome is a metric on $\O_Y(-P)$
whose curvature form restricts to $\o_P$ on $P$. Direct computation
shows that for $d \gg 0$, the curvature of the resulting tensor
product metric $||\cdot||_F$ on $F$ is a positive $(1,1)$-form,
which we denote by $\o_Y$.

For our application, we suppose that $||\cdot||_E$ has been chosen
in such a way that the two branches of $S$ meet orthogonally at $p$.
This is always possible, in fact one can modify the K{\"a}hler
potential to make the metric standard in any given local holomorphic
coordinates; see e.g.\ \cite[Lemma 7.2]{ruan02b}. The advantage is
that we can then identify $P \iso \P^2$ in such a way that
$C_0,C_1,C_2$ become the coordinate lines, {\em and} $\o_P$ the
standard Fubini-Study form. As a consequence, any Clifford torus
$K_0 \subset P$ is Lagrangian for $\o_Y|P = \o_P$. The K{\"a}hler
form $\o_Y$ also induces a symplectic connection on the set of
$\pi$-regular points in $Y$. We use this connection to transport
$K_0$ into nearby fibres, as described above. The resulting $K$ is a
Lagrangian submanifold with boundary inside $(Y,\o_Y)$ (this is best
seen in two steps: since the connection is symplectic and $K_0$ is
Lagrangian, each $K_t$ will be Lagrangian in $Y_t$; and since the
horizontal subspace is defined as the $\o_Y$-orthogonal complement
to the fibrewise tangent spaces, $K$ itself is Lagrangian).

\begin{lemma} \label{th:gromov}
Suppose that there is a sequence $t_k \in (0;1]$ with $\lim_k t_k =
0$, and a sequence of holomorphic discs $u_k: (\D,\partial \D)
\rightarrow (Y_{t_k},K_{t_k})$ whose areas $\int u_k^*\o_Y$ are
bounded. Then, after passing to a subsequence, there is a finite
collection of holomorphic maps $v_i: \P^1 \rightarrow X$ with the
property that (with respect to the isomorphism $Y_{t_k} \iso X$) the
image of $u_k$ for $k \gg 0$ is contained in an arbitrarily small
neighbourhood of the union of the images of the $v_i$. Moreover, the
union of the images of the $v_i$ is connected, and contains $p$.
\end{lemma}

\proof Consider $(t_k,u_k)$ as a sequence of holomorphic discs in
$Y$ with boundary on the Lagrangian submanifold $K$, and apply
Gromov compactness to a suitable subsequence. Since the images of
the discs lie in $Y_{t_k}$, the limiting stable disc has image in
$Y_0$. Its components are of three kinds: holomorphic spheres $w_i:
\P^1 \rightarrow Z$ and $y_i: \P^1 \rightarrow P$, as well as discs
$z_i: (\D,\partial \D) \rightarrow (P,K_0)$. A fairly weak
implication of Gromov convergence is that the image of $(t_k,u_k)$
for $k \gg 0$ is contained in an arbitrarily small neighbourhood (in
$Y$) of the union of the images of the $w_i,y_i,z_i$. We define the
$v_i$ to be the images of all the original components under blowdown
$Y \rightarrow \P^1 \times X$; the $y_i$ and $z_i$ become constant,
and in the latter case we replace the domain $\D$ by $\P^1$. The
convergence statement then holds by construction; since the original
stable disc was connected, the same applies to the union of the
images of the $v_i$; and since there was at least one $z_i$
component, there is at least one $v_i$ which is the constant map
with value $p$. \qed

\begin{lemma} \label{th:nef}
Suppose that there is a union of irreducible components of $S$,
forming a sub-curve $S' \subset S$ with $p \notin S'$, and an
effective {\em nef} divisor $D'$ whose support is $S'$. Assume that
we have $(t_k,u_k)$ as in the previous Lemma, with the additional
assumption that $u_k^{-1}(S') = \emptyset$ for all $k$. Then
$v_i^{-1}(S') = \emptyset$ for all $i$.
\end{lemma}

\proof We have $u_k \cdot D' = 0$ because the supports are disjoint,
and by looking at the Gromov limiting process, $\sum_i v_i \cdot D'
= 0$ (the fact that $(0,p)$ is blown up in $Y$ plays no role here,
since $p \notin S'$). Nefness implies that $v_i \cdot D' = 0$ for
each $i$, which means that the image of $v_i$ is either contained in
$S'$ or disjoint from it. Connectedness of the Gromov limit,
together with the fact that $p \notin S'$ lies on one of the $v_i$,
means that the first possibility is excluded. \qed

We now link this with the previous discussion of Stein-essential
Lagrangian submanifolds. $E$ comes with a canonical holomorphic
section $s_E$ which vanishes precisely on $S$. From this, the
section of $\O_{\P^1}(1)$ which vanishes exactly at $\{\infty\}$,
and the nowhere zero meromorphic section of $\O_Y(-P)$ which has a
simple pole at $P$, one gets a section $s_F$ of $F$ which vanishes
precisely on $Y_\infty \cup T \cup P$. For $t \neq 0,\infty$, choose
the isomorphism $F|Y_t \iso E$ in such a way that $s_F$ gets mapped
to $s_E$. Then the restriction of $||\cdot||_F$ to $Y_t \iso X$
induces a metric $||\cdot||_{E,t}$ on $E$, hence an exhausting
plurisubharmonic function $\phi_{M,t} = -\log ||s_E||_{E,t}$, which
by Lemma \ref{th:normal-crossing} makes $M$ into a finite type Stein
manifold. Note that the associated convex symplectic structure
$(\o_{M,t},\theta_{M,t})$ satisfies $\theta_{M,t} =
(d^c\log||s_F||_F)\, |\, (Y_t \setminus T_t)$ and $\o_{M,t} =
\o_Y|\,(Y_t \setminus T_t)$, where we are again using the
identifications $Y_t \setminus T_t \iso X \setminus S = M$.

By Lemma \ref{th:linking-torus} one can find a $\tau>0$ such that
for $t \in (0;\tau] \subset \P^1$, $K_t$ is a linking torus, and in
particular disjoint from $S$. This yields a smooth family of
submanifolds $K_t \subset M$ which are $\o_{M,t}$-Lagrangian.
Moreover, the class $[\theta_{M,t}|K_t]$ is constant in $t$. To see
this, note that since $K$ itself is $\o_Y$-Lagrangian, the
restriction of $d^c\log ||s_F||_F$ to $K \cap \pi^{-1}(0;\tau]$ is a
closed one-form. Its image under the restriction map $H^1(K \cap
\pi^{-1}((0;\tau]);\R) \rightarrow H^1(K_t;\R)$ is
$[\theta_{M,t}|K_t]$, which is therefore independent of $t$ as
claimed.

\begin{lemma} \label{th:non-essential}
Suppose that $K_\tau$, as a Lagrangian submanifold of the Stein
manifold $(M,\phi_{M,\tau})$, is {\em not} Stein-essential in the
sense of Section \ref{sec:essential}. Then there is a sequence $t_k
\in (0;\tau]$ with $\lim_k t_k = 0$, and a sequence of non-constant
holomorphic discs $u_k: (\D,\partial \D) \rightarrow (Y_{t_k}
\setminus T_{t_k},K_{t_k})$ whose areas $\int u_k^*\o_Y$ are
bounded.
\end{lemma}

\proof Suppose that the conclusion is false. Then as $t \in
(0;\tau]$ goes to zero, the least area of non-constant holomorphic
discs bounding $K_t$ in $(M,\o_{M,t})$ must go to infinity. More
precisely, for each $E>0$ there is a $\tau' \in (0;\tau]$ such that
every holomorphic disc $(\D,\partial \D) \rightarrow (M,K_{\tau'})$
with $\int u^* \o_{M,\tau'} \leq E$ is constant. The manifold $M$,
with its given complex structure and the family of plurisubharmonic
functions $\phi_{M,t}$, $t \in [\tau';\tau]$, is a finite type Stein
deformation (see the remark following Lemma
\ref{th:normal-crossing}), and the $K_t$, $t \in [\tau';\tau]$, are
a family of Lagrangian submanifolds such that $[\theta_{M,t}|K_t]$
is constant. By definition, the existence of such a deformation for
each $E$ means that $K_\tau$ is Stein-essential, contrary to our
assumption. \qed

Slightly more generally, suppose that for some $m \geq 1$, the
product $K_\tau^m$ is not Stein-essential as a Lagrangian
submanifold of $M^m$ equipped with the product Stein structure
(meaning the product complex structure and the plurisubharmonic
function $(x_1,\dots,x_m) \mapsto \phi_{M,\tau}(x_1) + \cdots +
\phi_{M,\tau}(x_m)$, which is still of finite type). Then the same
argument as before shows that one can find $t_k$ and non-constant
holomorphic discs $(u_k^1,\dots,u_k^m): (\D,\partial \D) \rightarrow
((Y_{t_k} \setminus T_{t_k})^m,K_t^m)$ with bounded area. After
choosing a non-constant component $u_k = u_k^{i_k}$, $i_k \in
\{1,\dots,m\}$ of each disc, one arrives at the same conclusion as
in the Lemma itself.

\section{Conclusion}\label{sec:conclusion}

We briefly recall Ramanujam's construction \cite{ramanujam71}. In
$\P^2$ take a smooth conic, and a cubic with a cusp singularity,
which intersect each other at two points with multiplicities $1$ and
$5$ respectively (the intersection points should also be distinct
from the cusp). Blow up the multiplicity $1$ intersection point, and
let $\bar{S}'$, $\bar{S}''$ be the proper transforms of the conic
and the cubic, respectively, inside the blowup $\bar{X} \iso \F_1$.
The union $\bar{S} = \bar{S}' \cup \bar{S}''$ has two singular
points, namely the multiplicity $5$ intersection point and the cusp.
Take the minimal blowup $q: X \rightarrow \bar{X}$ such that $S =
q^{-1}(\bar{S})$ is a divisor with normal crossing. The resolution
graph describing $S$ is shown in Figure \ref{fig:graph}. The
fattened vertices correspond to those components which are the
proper transforms of $\bar{S}'$ (with selfintersection $-2$) and
$\bar{S}''$ (with selfintersection $-3$); the other components are
exceptional divisors lying above the cusp (on the left) and the
multiplicity $5$ intersection point (on the right). Let $p \in S$ be
the crossing point in the preimage of the cusp where the proper
transform of the cubic intersects the exceptional divisor of
multiplicity $6$; this corresponds to the edge indicated by the
arrow.
\begin{figure}
\begin{center}
\begin{picture}(0,0)%
\includegraphics{graph.pstex}%
\end{picture}%
\setlength{\unitlength}{3947sp}%
\begingroup\makeatletter\ifx\SetFigFont\undefined%
\gdef\SetFigFont#1#2#3#4#5{%
  \reset@font\fontsize{#1}{#2pt}%
  \fontfamily{#3}\fontseries{#4}\fontshape{#5}%
  \selectfont}%
\fi\endgroup%
\begin{picture}(4353,1272)(526,-898)
\put(4726,-736){\makebox(0,0)[lb]{\smash{{\SetFigFont{12}{14.4}{\rmdefault}{\mddefault}{\updefault}{\color[rgb]{0,0,0}-2}%
}}}}
\put(526,-736){\makebox(0,0)[lb]{\smash{{\SetFigFont{12}{14.4}{\rmdefault}{\mddefault}{\updefault}{\color[rgb]{0,0,0}-3}%
}}}}
\put(1126,-736){\makebox(0,0)[lb]{\smash{{\SetFigFont{12}{14.4}{\rmdefault}{\mddefault}{\updefault}{\color[rgb]{0,0,0}-1}%
}}}}
\put(1726,-736){\makebox(0,0)[lb]{\smash{{\SetFigFont{12}{14.4}{\rmdefault}{\mddefault}{\updefault}{\color[rgb]{0,0,0}-3}%
}}}}
\put(2326,-736){\makebox(0,0)[lb]{\smash{{\SetFigFont{12}{14.4}{\rmdefault}{\mddefault}{\updefault}{\color[rgb]{0,0,0}-1}%
}}}}
\put(2326,239){\makebox(0,0)[lb]{\smash{{\SetFigFont{12}{14.4}{\rmdefault}{\mddefault}{\updefault}{\color[rgb]{0,0,0}-2}%
}}}}
\put(1126,239){\makebox(0,0)[lb]{\smash{{\SetFigFont{12}{14.4}{\rmdefault}{\mddefault}{\updefault}{\color[rgb]{0,0,0}-2}%
}}}}
\put(2926,-736){\makebox(0,0)[lb]{\smash{{\SetFigFont{12}{14.4}{\rmdefault}{\mddefault}{\updefault}{\color[rgb]{0,0,0}-2}%
}}}}
\put(3526,-736){\makebox(0,0)[lb]{\smash{{\SetFigFont{12}{14.4}{\rmdefault}{\mddefault}{\updefault}{\color[rgb]{0,0,0}-2}%
}}}}
\put(4126,-736){\makebox(0,0)[lb]{\smash{{\SetFigFont{12}{14.4}{\rmdefault}{\mddefault}{\updefault}{\color[rgb]{0,0,0}-2}%
}}}}
\end{picture}%
\caption{\label{fig:graph}}
\end{center}
\end{figure}

Because $\bar{S}$ is ample, one can find an ample divisor $D$ on $X$
whose support is $S$. Set $E = \O_X(D)$, and carry out the
construction from the previous section for the point $p$; this
yields a family of plurisubharmonic functions $\phi_{M,t}$ and
$\o_{M,t}$-Lagrangian tori $K_t$. Take $\tau$ as in the discussion
before Lemma \ref{th:non-essential}. Equip $M$ with its standard
complex structure $J_M$, the function $\phi_M = \phi_{M,\tau}$ which
makes it into a finite type Stein manifold, and the Lagrangian
submanifold $L = K_\tau$.

Let $S' \subset S$ be the preimage of $\bar{S}'$, which is its
proper transform together with the exceptional curves arising from
blowing up the multiplicity $5$ intersection point. Since $\bar{S}'
\subset \bar{X}$ is ample, one can use Kodaira's construction to
find an effective divisor $D'$ with support $S'$, and a $(1,1)$-form
representing its Poincar\'e dual, which is nonnegative everywhere,
and positive away from the preimage of the cusp point. This means
that any curve $\Sigma$ with $\Sigma \cdot D' \leq 0$ must be lie on
the preimage of the cusp point; in particular $D'$ is nef.

\proof[Proof of Theorem \ref{th:main}] Assume that $M^m$ has a
symplectic embedding $i$ into a subcritical Stein manifold $N$,
which we may assume to be complete by Lemma \ref{th:stupid}. Since
$H^1(M^m;\R) = 0$, the closed one-form $i^*\theta_N - \theta_M$ is
automatically exact. By Proposition \ref{th:proposition}, $L^m$
cannot be Stein-essential. Using Lemma \ref{th:non-essential} and
the remark following it, one then gets a sequence $t_k$ and
non-constant holomorphic discs $u_k: (\D,\partial \D) \rightarrow
(Y_{t_k} \setminus T_{t_k},K_{t_k})$ with bounded energy. The
isomorphism $Y_{t_k} \iso X$ sends $T_{t_k}$ to $S$, hence the image
of each $u_k$ is disjoint from $S' \subset S$. The limit of this, in
the sense of Lemma \ref{th:gromov}, is a finite collection of
holomorphic maps $v_i: \P^1 \rightarrow X$, which by Lemma
\ref{th:nef} are disjoint from $S'$. Using the observation made
above, it follows that the image of the $v_i$ lies in the preimage
of the cusp point. In other words, if $\bar{V}$ is a neighbourhood
of the cusp in $\bar{X}$, and $V$ its preimage in $X$, then for $k
\gg 0$ we have a non-constant holomorphic disc $u_k: (\D,\partial
\D) \rightarrow (V \setminus S,K_{t_k})$. Since $K_{t_k}$ is a
linking torus for a crossing point which arose from a cusp, we can
apply Lemma \ref{th:trefoil} to conclude that $u(\partial \D)$ is a
contractible loop on $K_{t_k}$. But by Stokes this implies that the
area $\int u_k^*\o_Y$ is zero, which means that $u_k$ is constant,
hence yields a contradiction. \qed

\proof[Proof of Corollary \ref{th:cor}] Suppose that $M^m$ carries
an exhausting plurisubharmonic Morse function $\tilde\phi_M$ with
just $1$ critical point. In view of Lemma \ref{th:stein-rescale}, we
may assume that $\tilde\phi_M$ is complete; and it is of finite type
by assumption. Lemma \ref{th:two-functions} then says that there is
a symplectic embedding $i: M \rightarrow M$ such that
$i^*\tilde\theta_M - \theta_M$ is exact, contradicting Theorem
\ref{th:main}. \qed

{\smaller
\section*{Appendix}

\proof[Proof of Lemma \ref{th:weinstein-embed}] This is a
straightforward generalization of the case of cotangent bundles,
discussed in \cite[Proposition 1.3.A]{eliashberg-gromov93}. Without
affecting the validity of the statement, we may replace $\theta_N$
by $\tilde{\theta}_N = \theta_N + dK$ for any compactly supported
function $K$. A suitable choice of $K$ ensures that
$i^*\tilde\theta_N = \theta_M$, and then $i$ takes the Liouville
flow $l_M$ to the modified Liouville flow $\tilde{l}_N$ associated
to $\tilde{\theta}_N$. By assumption, for any point $x \in M$ there
is a $t \geq 0$ such that $l_M^{-t}(x) \in
\phi_M^{-1}((-\infty;c_0])$, and on the other hand $\tilde{l}_N^t$
is defined for all $t \geq 0$. Hence \[j_t = \tilde{l}_N^t \circ i
\circ l_M^{-t}, \qquad t \geq 0\] is a family of mutually compatible
extensions of $i$ to successively larger subsets, which exhaust $M$;
and they satisfy $j_t^*\tilde{\theta}_N = \theta_M$. \qed

\proof[Proof of Lemma \ref{th:weinstein-deform}] After a finite
decomposition of the interval $[0;1]$, we may assume that there are
$c_1<c_2<\dots$ converging to $+\infty$, such that
$d\phi_{M,t}(\lambda_{M,t}) > 0$ on $\phi_{M,t}^{-1}(c_k)$ for all
$t \in [0;1]$. There is a $k$ such that $U \subset
\phi_{M,0}^{-1}((-\infty;c_k])$. Let $\mu_t$ be the Moser vector
field defined by $\o_{M,t}(\mu_t,\cdot) = -d\theta_{M,t}/dt$. By
choosing $r \gg 0$ sufficiently large, one can achieve that $
d\phi_{M,t}(\mu_t - r\lambda_{M,t}) + \partial_t \phi_{M,t} < 0 $ on
$\phi_{M,t}^{-1}(c_k)$ for all $t \in [0;1]$. Hence, integrating
$\mu_t - r\lambda_{M,t}$ yields a smooth family of embeddings $i_t:
\phi_{M,0}^{-1}((-\infty;c_k]) \rightarrow
\phi_{M,t}^{-1}((-\infty;c_k])$ starting with $i_0 = id$, such that
$i_t^*\theta_{M,t} - e^{-rt}\theta_{M,0}$ is an exact one-form.
Define $j_t$ by composing $i_t|U$ and the time $rt$ flow of
$\lambda_{M,t}$. \qed

\proof[Proof of Lemma \ref{th:weinstein-deform-2}] Let $c_0$ be as
in the definition of finite type deformations.  Gray's Theorem on
the stability of contact structures implies that there is a family
of diffeomorphisms of $\phi_{M,t}^{-1}(c_0)$ which pulls back (the
restrictions of) $\theta_{M,t}$ to $\theta_{M,0}$. This induces a
family of diffeomorphisms which identify the cone-like ends of
$(M,\o_{M,t},\theta_{M,t})$ for different $t$. Going back from
isotopies to the generating vector fields, the outcome is that one
can find vector fields $\gamma_t$ on $\phi_{M,t}^{-1}([c_0;\infty))$
such that $L_{\gamma_t}\theta_{M,t} + \partial_t\theta_{M,t} = 0$
and $[\gamma_t,\lambda_{M,t}] + \partial_t\lambda_{M,t} = 0$ (the
first condition comes from Gray's argument, and the second one
expresses compatibility with the conical structure). Extend these
vector fields to all of $M$ in an arbitrary way, and integrate them
to get a family of diffeomorphisms $g_t: M \rightarrow M$ such that
$g_t^*\theta_{M,t} - \theta_{M,0}$ vanishes outside a compact
subset. Moser's Lemma then yields another family of diffeomorphisms
$h_t: M \rightarrow M$, which is compactly supported and satisfies
$h_t^* g_t^*\theta_{M,t} - \theta_{M,0} = dR_t$ as desired. Set $f_t
= g_t \circ h_t$. \qed

\proof[Proof of Lemma \ref{th:stein-rescale}] The main statement is
taken from \cite[Lemma 3.1]{biran-cieliebak01}. There, the authors
observe that for any $h$ with $h' > 0$, $h'' \geq 0$, the modified
function $\tilde\phi_M = h(\phi_M)$ is again plurisubharmonic.
Moreover, the Liouville vector field associated to the modified
Stein structure is related to the original one by
\[
 \tilde\lambda_M =
 \lambda_M \cdot \frac{h'(\phi_M)}{h'(\phi_M) + h''(\phi_M)
 ||\nabla\phi_M||^2},
\]
with the norm taken in the original K{\"a}hler metric. This means
that the modified Liouville flow has the same flow lines as the
original one, but moves along them at a slower or equal rate. The
additional condition $h''(c) \geq \delta h'(c)$ for $c \geq c_0$
implies $d\phi_M(\tilde\lambda_M) \leq \delta^{-1}$ outside a
compact subset, so that the flow is then complete.

To prove the last statement in the Lemma, consider the family of
functions $\phi_{M,t} = h_t(\phi_M)$ with $h_t(c) = (1-t)c + th(c)$,
$0 \leq t \leq 1$. These also satisfy $h_t' > 0$, $h_t'' \geq 0$, so
the flow of $\nabla \phi_{M,t}$ (with respect to its associated
K{\"a}hler metric) is slower than that of $\nabla \phi_M$. We are
assuming the second flow is complete, hence so is the first one.
Besides that, we are also assuming that $\phi_M$ is of finite type.
Hence what we get is a complete finite Stein deformation
$(M,J_M,\phi_{M,t})$, to which Lemma \ref{th:weinstein-deform-2} can
be applied. \qed

\proof[Proof of Lemma \ref{th:two-functions}] We imitate the
argument from \cite[Theorem 1.4.A]{eliashberg-gromov93}, with some
clarifications. We will prove the statement first in the case where
$\tilde\phi_M$ grows faster than $\phi_M$, by which we mean that the
difference $\delta = \tilde\phi_M - \phi_M$ is an exhausting
function. Let $U_k = \delta^{-1}((-\infty;k))$ be the associated
family of exhausting subsets; after changing $\tilde\phi_M$ by a
constant, we may assume that $\delta \geq 1$ everywhere, hence $U_k
= \emptyset$ for $k \leq 1$. Fix a smooth function $l: \R
\rightarrow \R$ such that $l(r) = 0$ for $r \leq -1$, $l(r) = r$ for
$r \geq 1$, and $l''(r) \geq 0$ everywhere. For $k = 0,1,2,\dots$
consider the functions $\tilde\phi_{M,k} = \phi_M + l(\tilde\phi_M -
\phi_M - k)$. These satisfy $\tilde\phi_{M,k} = \phi_M$ on
$\bar{U}_{k-1}$, $\tilde\phi_{M,k} = \tilde\phi_M - k$ on $M
\setminus U_{k+1}$, and are plurisubharmonic because
\[
 dd^c \tilde\phi_{M,k} = (1-l')\, dd^c\phi_M + l'
 \, dd^c\tilde{\phi}_M
 + l'' \, d(\tilde\phi_M - \phi_M)
 \wedge d^c(\tilde\phi_M - \phi_M)
\]
with $l' = l'(\tilde\phi_M - \phi_M - k) \in [0;1]$, and $l'' =
l''(\tilde\phi_M - \phi_M - k) \geq 0$. Let
$(\tilde\o_{M,k},\tilde\theta_{M,k})$ be the convex symplectic
structure associated to $\tilde\phi_{M,k}$. By applying Moser's
argument to the linear deformation between the $k^{\text{th}}$ and
$(k+1)^{\text{st}}$ of these structures, we get a a diffeomorphism
$f_k: M \rightarrow M$ which is the identity outside $\bar{U}_{k+2}
\setminus U_{k-1}$, and such that $f_k^*\tilde\theta_{M,k} -
\tilde\theta_{M,k+1}$ is the derivative of a function $K_k$
supported in $\bar{U}_{k+2} \setminus U_{k-1}$. Let $i: M
\rightarrow M$ be the infinite composition $f_0 \circ f_1 \circ
\cdots$. This is well-defined because for each $x \in M$, one has
$f_k(x) = x$ for all but finitely many $k$. The infinite composition
is injective and a local diffeomorphism, hence an embedding (but not
necessarily a diffeomorphism; composing the $f_k^{-1}$ in the
opposite order makes no sense). By definition $\tilde\theta_{M,0} =
\tilde\theta_M$; and for each relatively compact subset $U \subset
M$ there is a $k$ such that $\tilde\theta_{M,k} = \theta_M$ on $U$,
and $f_{k+1}|U = f_{k+2}|U = \cdots = id_U$. It follows that
$i^*\tilde\o_M = \o_M$, and moreover $i^*\tilde\theta_M = \theta_M +
dK$ for $K = (f_1 \circ f_2 \circ \cdots)^* K_0 + (f_2 \circ f_3
\circ \cdots)^* K_1 + \cdots$, the same argument as before showing
that this sum is well-defined.

We now pass to the general situation, where $\tilde\phi_M$ is
complete and of finite type but otherwise arbitrary. One can then
find a function $h$ as in Lemma \ref{th:stein-rescale} for which
$h(\tilde\phi_M)$ grows faster than $\phi_M$, and moreover this
rescaling does not change the (exact) symplectic isomorphism type.
One can then apply the previous argument to $\phi_M$ and
$h(\tilde\phi_M)$, and derive the desired result. \qed

\proof[Proof of Lemma \ref{th:normal-crossing}] For simplicity, we
write $||\cdot|| = ||\cdot||_E$ and $s = s_E$. Around a point $x \in
S$, choose local holomorphic coordinates, and a local holomorphic
trivialization of our line bundle, with respect to which $s(z) =
z_1^{w_1} \cdots z_n^{w_n}$. Write $w = w_1 + \cdots + w_n$. With
respect to the trivial metric $||\cdot||_0$ one has
\[
\begin{aligned}
 |d||s||_0| & \; \wigglygeq \sum_j w_j \left| \partial_{z_j} ||s||_0 \right|
 \wigglygeq \sum_j w_j |z_1^{w_1} \cdots z_j^{w_j-1} \cdots
 z_n^{w_n}| \\
 & \; \wigglygeq \sqrt[w]{|z_1|^{w_1(w-1)} \cdots
 |z_n|^{w_n(w-1)}}
 = ||s||_0^{1-1/w},
\end{aligned}
\]
where $\wigglygeq$ means greater or equal than some small constant
times the right hand side (in spite of that, we have kept the $w_j$,
because they indicate how the inequality between arithmetic and
geometric mean is applied). One also has $|d ||s||\, | + ||s||
\wigglygeq |d ||s||_0|$ and $||s||_0^{1-1/w} \wigglygeq
||s||^{1-1/w}$. After combining the inequalities, one sees that $d
||s||$ does not vanish at points $z$ where $||s(z)||$ is
sufficiently small, hence $x$ does not lie in the closure of the
critical point set of $\phi_M$. \qed

} % end of "smaller"

%\bibliographystyle{amsplain}
%\bibliography{../../../bib/all,../../../bib/new}
\providecommand{\bysame}{\leavevmode\hbox
to3em{\hrulefill}\thinspace}
\providecommand{\MR}{\relax\ifhmode\unskip\space\fi MR }
\providecommand{\MRhref}[2]{%
  \href{http://www.ams.org/mathscinet-getitem?mr=#1}{#2}
} \providecommand{\href}[2]{#2}

\end{document}